\documentclass[12pt, reqno]{amsart}
\usepackage{amsmath, amsthm, amscd, amsfonts, amssymb, graphicx, color}
\usepackage[bookmarksnumbered, plainpages]{hyperref}
\usepackage{mathtools}
\usepackage{scalerel}
\usepackage{nccmath}
\usepackage{bm}
\usepackage{commath}

\newtheoremstyle{break}
{\topsep}{\topsep}%
{\itshape}{}%
{\bfseries}{}%
{\newline}{}%

\theoremstyle{break}
\newtheorem{definition}{Definition}[section]
\newtheorem{theorem}[definition]{Theorem}
\newtheorem{lemma}[definition]{Lemma}

\newtheorem{corollary}[definition]{Corollary}

\newtheoremstyle{plain}
{\topsep}{\topsep}%
{}{}%
{\bfseries}{}%
{\newline}{}%

\theoremstyle{plain}

\newenvironment{ProofOfNo}{$ $}{$ $\null\hfill$\square$}

\newcommand{\D}{\mathbb{D}}
\newcommand{\C}{\mathbb{C}}
\newcommand{\F}{\mathcal{F}}
\newcommand{\N}{\mathbb{N}}
\newcommand{\R}{\mathbb{R}}

\newcommand{\fn}{(f_n)_n}
\newcommand{\HD}{\mathcal{H}(\D)}
\newcommand{\MD}{\mathcal{M}(\D)}

\begin{document}
	\setcounter{page}{1}
	\mathtoolsset{showonlyrefs,showmanualtags}
	
	\title{Lifts of Logarithmic Derivatives}	
	\author{Matthias Grätsch}
	\maketitle	
	\begin{abstract}
		Consider a sequence of meromorphic functions $\fn$. This paper presents a technique that enables the transfer of convergence properties from $(f_n^{(m+1)}/f_n^{(m)})_n$ to subsequences of $(f_n^{(m)}/f_n^{(m-1)})_n$. As an application, we will show that the families of functions with bounded Schwarzian derivative are quasi-normal. 
	\end{abstract} 
	\section{Introduction and Statement of Results}
	\noindent 
	Throughout this paper, we are going to denote the set of all holomorphic functions on a domain $D\subseteq\C$ as $\mathcal{H}(D)$. Likewise, we will write $\mathcal{M}(D)$ for the set of all meromorphic functions on $D$. Moreover, we denote the open unit disk by ${\D\coloneq\{z\in\C\,\,:\,\,\abs{z}<1\}}$ and the set of all poles of $f\in\mathcal{M}(D)$ by $\mathcal{P}_f$. Similarly, we denote the set of all zeros of $f$ by $\mathcal{Z}_f$. Whenever we are counting the multiplicity of such a zero or pole, we will remark ``(CM)'' and we will remark ``(IM)'' if we ignore it.\medskip
	
	\noindent
	We begin with a simple example to illustrate the main purpose of this paper. Consider the sequence of functions $f_n(z)=nz^2+1$. Then we get for the following sequences that:
	\begin{align}
		&\frac{f_n'''(z)}{f_n''(z)}=0&&\text{ converges with an analytic limit.}\\
		&\frac{f_n''(z)}{f_n'(z)}=\frac{1}{z}&&\text{ converges with a meromorphic limit.}\\
		&\frac{f_n'(z)}{f_n(z)}=\frac{2nz}{nz^2+1}&&\!
		\begin{array}{l}
			\text{does not converge locally uniformly at }0 \text{ but}\\
			\text{it converges on }\C\backslash\{0\}\text{ with an analytic limit}.
		\end{array}
	\end{align}
	\noindent
	Note that all three sequences converge on the largest part of their domains. The basic idea of this paper is that the convergence of $(f_n^{(m+1)}/f_n^{(m)})_n$ to an analytic limit induces (to some extent) the convergence of a subsequence of $(f_n^{(j+1)}/f_n^{(j)})_n$ for ${0\leq j<m}$. However, the ``quality of convergence'' gradually deteriorates for logarithmic derivatives of lower order, which allows for poles in the limit and some points where the logarithmic derivatives of lower order do not converge at all.
	To quantify this decline, we need the concept of ${Q_m\text{-normality}}$, which was introduced by C.-T. Chuang in \cite{Chuang}.\pagebreak
	\begin{definition}
		Let $D\subseteq\C$ be a domain and $E\subseteq D$. The \textbf{derived set} of E with respect to D is the set of all accumulation points of $E$ in $D$. We denote this set as $E_D^{(1)}$. Inductively, we can define the \textbf{derived set of order k} $\in\N$ of $E$ with respect to $D$ as $E_D^{(k)}\coloneq(E_D^{(k-1)})_D^{(1)}$. We define $E_D^{(0)}\coloneq E$ as well.\\
		Let $m\in\N_0$. We say that a family $\F\subseteq\mathcal{M}(D)$ is \textbf{$\bm{Q_m}$-normal} on $D$ if for each sequence $(f_n)_n\subseteq\F$, there exists a set $E\subseteq D$ with $E_D^{(m)}=\emptyset$ and a subsequence $(f_{n_k})_{k}\subseteq\fn$ that converges locally uniformly (with respect to the spherical metric) on $D\backslash E$. 
		If $m\in\N_0$ is the smallest number such that $\F$ is $Q_m$-normal, then we will call $m$ the \textbf{index of normality} of $\F$.\\
		If a family is $Q_1$-normal, it is called \textbf{quasi-normal}.
	\end{definition} \medskip
	
	\noindent
	In our example, we get that $(f_n'''/f_n'')_n$ is $Q_0$-normal (i.e. normal) with an analytic limit. In that case, we are able to ``lift this property'' and get that $(f_n''/f_n')_n$ is $Q_0$-normal as well, but with the limit function having a pole. However, when we are lifting again to the sequence $(f_n'/f_n)_n$, the index of normality increases by one. In Lemma \ref{le:LogarithmicDerivativesZalcman}, we will see that this is the maximum increase in the index of normality that can occur when we lift from one logarithmic derivative to another. Thus, in the general case, we get that an analytic or meromorphic limit of $(f_n^{(m+1)}/f_n^{(m)})_n$ implies for all $j=0,\hdots, m-1$ that a subsequence of $(f_n^{(j+1)}/f_n^{(j)})_n$ converges almost everywhere in its domain due to the countability of the exceptional sets. Lemma \ref{le:LogarithmicDerivativesZalcman} will also show that the limits of these subsequences are always meromorphic and that none of these subsequences converge to $\infty$. \medskip
	
	\begin{theorem}\label{th:Lift}
		Let $D\subseteq\C$ be a domain, $m\in\N$ and $\fn\subseteq\mathcal{M}(D)$, such that $(f_n^{(m+1)}/f_n^{(m)})_n$ converges locally uniformly to some $F_m\in\mathcal{M}(D)$. Then we get for $0\leq j\leq m$ that $(f_n^{(j+1)}/f_n^{(j)})_n$ is $Q_{(m-j)}$-normal.\\
		In particular, there exists a subsequence $(f_{n_k})_k\subseteq\fn$ as well as an open and dense subset $U\subseteq D$ where $(f_{n_k}^{(j+1)}/f_{n_k}^{(j)})_k$ converges locally uniformly to an analytic limit for $j=0,\hdots,m$.
	\end{theorem}\medskip
	
	\noindent
	The example shows that the index of normality does not have to increase when we are lifting. In fact, it is unknown to the author if lifting can produce an index of normality $\geq\!2$, i.e if there is $m\in\N$ and a sequence $(g_n)_n\subseteq\MD$, where $(g_n^{(m+1)}/g_n^{(m)})_n$ converges locally uniformly to some $G_m\in\MD$ while $(g_n'/g_n)$ is $Q_2$-normal and not quasi-normal.\\
	\noindent
	Note that the convergence of $(g_n^{(m+1)}/g_n^{(m)})_n$ to some $G_m\in\MD$ implies that $(g_n)_n$ is $Q_2$-normal. This is  due to a recent paper by S. Nevo and Z. Shem Tov. They have shown for any domain $D$ that a family $\F\subseteq\mathcal{H}(D)$ that satisfies 
	\begin{align}\label{eq:NevoTov}
		\frac{\mathopen|f^{(m+1)}\mathclose|}{1+\mathopen|f^{(m)}\mathclose|}(z)\leq M \qquad\text{for all } f\in \F\text{ and }z\in D
	\end{align}
	for some $m\in\N$ and $M\in\R$ is already quasi-normal (see \cite[{Theorem 2}]{Nevo}).
	Note that $(g_n)_n$ fulfills \eqref{eq:NevoTov} locally for all points $z\in\D\backslash\mathcal{P}_{G_m}$ for almost all $n\in\N$. Since ${Q_m\text{-normality}}$ is a local property (see \cite[Theorem 8.3]{Chuang}), we are able to conclude that $(g_n)_n$ is $Q_1$-normal on $\D\backslash\mathcal{P}_{G_m}$ and $Q_2$-normal on $\D$. This is why, in the author's opinion, it is unlikely that higher indices of normality occur during lifting. To illustrate this point further, this paper will also present several additional conditions, in which the index of normality cannot increase during lifting.\\
	In the proof of Theorem \ref{th:Lift} (and more specifically Lemma \ref{le:LogarithmicDerivativesZalcman}), we will see that the exceptional set $E$ which corresponds to a convergent subsequence $(f_{n_k}^{(j+1)}/f_{n_k}^{(j)})_k$ depends mainly on the poles of $(f_{n_k}^{(j+2)}/f_{n_k}^{(j+1)})_k,\hdots,(f_{n_k}^{(m+1)}/f_{n_k}^{(m)})_k$ and their accumulation points. Note that all poles of $f_n$ are also poles of $f_n^{(m+1)}/f_n^{(m)}$, and the set consisting of them is small in the sense of set derivatives since $(\mathcal{P}_{F_m})_D^{(1)}=\emptyset$. So we can say that the index of normality depends mainly on the zeros of the derivatives $f_{n_k}^{(j)},\hdots,f_{n_k}^{(m)}$. With additional constraints on these zeros, it is possible to characterize some cases where the index of normality does not increase during lifting. One example is the following corollary:\bigskip
	\begin{corollary}\label{co:Extension}
		Let $D\subseteq\C$ be a domain, $N, m\in\N$ with $m\geq2$ and let $\fn\subseteq\mathcal{M}(D)$ be a sequence with:
		\begin{enumerate}
			\item[\textit{(a)}] $\Big(\mfrac{f_n^{(m+1)}}{f_n^{(m)}}\Big)_n$ converges locally uniformly to some $F_m\in\mathcal{M}(D)$ on $D$.
			\item[\textit{(b)}] The multiplicity of every zero of $f_n^{(m)}$ is at most $N$.
		\end{enumerate}
		Then we know that $(f_n^{(m)}/f_n^{(m-1)})_n$ and $(f_n^{(m-1)}/f_n^{(m-2)})_n$ are both quasi-normal and none of their subsequences converges to $\infty$.
	\end{corollary}	\bigskip
	\noindent
	As an application of such lifts of logarithmic derivatives, we are going to show the following result about functions with bounded Schwarzian derivative 
	\begin{align}\nonumber
		S_f\coloneq\bigg(\frac{f''}{f'}\bigg)'-\frac{1}{2} \bigg(\frac{f''}{f'}\bigg)^2.
	\end{align}
	\medskip
	\begin{theorem}\label{th:Schwarz}
		For any $M\in\R$, we get that $\F_M\coloneq\{f\in\MD\;:\,||S_f||_\infty\leq M\}$  is quasi-normal.
	\end{theorem}
	\bigskip\noindent
	The properties of analytic functions with bounded Schwarzian were already studied in \cite{MaMeijaMinda} by ${\text{W. Ma}}$, ${\text{D. Meija}}$ and D. Minda. Thus, in the special case of holomorphic functions, Theorem \ref{th:Schwarz} can be obtained from ${\text{\cite[Theorem 1.1(b)]{MaMeijaMinda}}}$ in combination with \cite[Corollary 5.9]{Chuang}.\\
	However, we are going to take a different approach. We will show that the family of pre-Schwarzians $\F''_M/\F'_M\coloneq\{f''/f'\,:\,f\in\F_M\}$ is normal and that none of its sequences converges to $\infty$. Then, we will lift this property to $(f_n'/f_n)_n$ and $\fn$. Here, we are going to use the well-known fact that $f'$ is zero-free if $S_f\in\HD$ to keep the index of normality low enough. An extensive collection of such results about the Schwarzian derivative can be found in \cite[Chapter II]{Letho} by O. Lehto. 	\newpage
	
	\section{Auxiliary Lemmas and their Proofs}
	\noindent
	All results in this paper depend heavily on Zalcman's Lemma.
	\medskip
	\begin{lemma}[\cite{Zalcman}, Zalcman's Lemma]\label{ZalcmanLemma}
		Let $\mathcal{F}\subseteq\MD$ be a family of meromorphic functions and suppose that $\F$ is not normal at $z_0\in\D$. Then there are sequences $\fn\subseteq\F$, $(z_n)_n\subseteq\D$ and $(\rho_n)_n\subseteq(0,1)$ with $\lim_{n\rightarrow\infty}\rho_n=0$ and $\lim_{n\rightarrow\infty} z_n=z_0$, such that 
		\begin{align*}
			g_n(\zeta)\coloneq f_n(z_n+\rho_n\zeta)
		\end{align*}
		converges locally uniformly to a non-constant $g\in\mathcal{M}(\C)$ with ${g^\#(\zeta)\leq g^\#(0)=1}$ for all $\zeta\in\C$.
	\end{lemma}
	\medskip\noindent
	The increase in the index of normality is due to the following lemma.
	\medskip
	\begin{lemma}\label{le:LogarithmicDerivativesZalcman}
		Let $D\subseteq\C$ be a domain, $m\in\N$ and $\fn\subseteq\mathcal{M}(D)$, such that $(f_n^{(m+1)}/f_n^{(m)})_n$ converges locally uniformly to $F_m\in\mathcal{M}(D)$. Then there exists a subsequence ${(f_{n_k})_{k}\subseteq\fn}$, such that $\big(f_{n_k}^{(m)}/f_{n_k}^{(m-1)}\big)_{k}$ converges locally uniformly on $D\backslash\mathcal{P}_{F_m}$ to a meromorphic function $F_{m-1}\in\mathcal{M}(D\backslash\mathcal{P}_{F_m})$.
	\end{lemma}
	\medskip\noindent	
	The usage of Zalcman's Lemma in this context is a bit reminiscent of \cite{GrahlNevo}. Some of the techniques we will use can already be found in there.\\
	Next, we introduce a new notion to abbreviate our later reasoning a little bit.
	\medskip
	\begin{definition}
		If $D\subseteq\C$ is a domain, we say that a sequence $\fn\subseteq\mathcal{M}(D)$ is \textbf{eventually bounded} on D if, for every point $w\in D$, there is a neighborhood $U_w\subseteq D$ around $w$ and there exist $N_w\in\N$ as well as $M_w\in\R$ such that $|f_n(z)|\leq M_w$ for all $z\in U_w$ and all $n\geq N_w$.\\
	\end{definition}
	\noindent
	\textit{Proof of Lemma \ref{le:LogarithmicDerivativesZalcman}:}
	\begin{ProofOfNo}
		Because of $\lim_{n\rightarrow\infty}f_n^{(m+1)}/f_n^{(m)}=F_m$, we know that the sequence $(f_n^{(m+1)}/f_n^{(m)})_n$ is eventually bounded on $D\backslash\mathcal{P}_{F_m}$. In other words, for every $w\in D\backslash\mathcal{P}_F$, there is a neighborhood $U_{w}$ around $w$ and $M_{w}\in\R$, such that for every $z\in U_{w}$, we have $\big|f_n^{(m+1)}(z)/f_n^{(m)}(z)\big|\leq M_w$ for almost all $n\in\N$.\\ \sloppy
		Next, we suppose that $(f_n^{(m)}/f_n^{(m-1)})_n$ is not normal in some point ${z_0\in D\backslash\mathcal{P}_{F_m}}$. Then, Zalcman's Lemma (Lemma \ref{ZalcmanLemma}) states that there are sequences ${(\rho_n)_n\subseteq(0,1)}$ and ${(z_n)_n\subseteq D\backslash\mathcal{P}_{F_m}}$, as well as a subsequence of $(f_n^{(m)}/f_n^{(m-1)})_n$ (that we will continue to call $(f_n^{(m)}/f_n^{(m-1)})_n$), such that
		\begin{align*}
			g_n(\zeta)\coloneq\frac{f_n^{(m)}}{f_n^{(m-1)}}(z_n+\rho_n\zeta)
		\end{align*}
		converges locally uniformly to a non-constant limit $g\in\mathcal{M}(\C)$.
		Now, we can consider a closed disk $B\subseteq\C$, such that $g$ has no pole in $B$. Note that there is $N\in\N$, such that $(g_n)_{n\geq N}$ is uniformly bounded on $B$, i.e. there is also a constant $C\in\R$ with ${|g_n(\zeta)|\leq C}$ for all $\zeta\in B$ whenever $n\geq N$. \\However, the Weierstraß convergence theorem states that for all $\zeta\in B$, we have
		\begin{align*}
			g_n'(\zeta)=\rho_n\Bigg(\frac{f_n^{(m+1)}}{f_n^{(m-1)}}-\bigg(\frac{f_n^{(m)}}{f_n^{(m-1)}}\bigg)^2\Bigg) (z_n+\rho_n\zeta)\rightarrow g'(\zeta)
		\end{align*}
		and respectively
		\begin{align*}
			g'(\zeta)&=\lim\limits_{n\rightarrow\infty}g_n'(\zeta)\\
			&=\lim\limits_{n\rightarrow\infty}\rho_n\Bigg(\underbrace{\frac{f_n^{(m)}}{f_n^{(m-1)}}(z_n+\rho_n\zeta)}_{=g_n(\zeta)}\frac{f_n^{(m+1)}}{f_n^{(m)}}(z_n+\rho_n\zeta)-\underbrace{\bigg(\frac{f_m^{(m)}}{f_n^{(m-1)}}(z_n+\rho_n\zeta)\bigg)^2}_{=g_n^2(\zeta)}\Bigg)\\
			&=\lim\limits_{n\rightarrow\infty}\underbrace{\rho_n\vphantom{\frac{f_n^{(m)}}{f_n^{(m)}}}}_{\rightarrow0} \underbrace{g_n(\zeta)\vphantom{\frac{f_n^{(m)}}{f_n^{(m)}}}}_{|\,\cdot\,|\,\leq\, C} \Bigg(\underbrace{\frac{f_n^{(m+1)}}{f_n^{(m)}}(z_n+\rho_n\zeta)}_{|\,\cdot\,|\,\leq\, M_{z_{\scaleto{0}{2.5pt}}}} - \underbrace{g_n(\zeta)\vphantom{\frac{f_n^{(m)}}{f_n^{(m)}}}}_{{|\,\cdot\,|\,\leq\, C}}\Bigg)
			=0.
		\end{align*}
		This contradicts Zalcman's Lemma. Therefore, we are able to conclude that $(f_n^{(m)}/f_n^{(m-1)})_n$ is in fact normal in $D\backslash\mathcal{P}_{F_m}$.\\
		Next, we  assume that there is a subsequence of $(f_n^{(m)}/f_n^{(m-1)})_n$ (that we are going to continue to call $(f_n^{(m)}/f_n^{(m-1)})_n$), which converges locally uniformly to $\infty$ on $D\backslash\mathcal{P}_{F_m}$. Then this would imply that both $(f_n^{(m-1)}/f_n^{(m)})_n$ and $\big((f_n^{(m-1)}/f_n^{(m)})'\big)_n$ have to converge locally uniformly to $0$ on $D\backslash\mathcal{P}_{F_m}$. However, this is a contradiction because for all $z\in\D\backslash\mathcal{P}_{F_m}$, we have
		\begin{align*}
			\lim\limits_{n\rightarrow\infty}\bigg(\frac{f_n^{(m-1)}}{f_n^{(m)}}\bigg)'(z) = \lim\limits_{n\rightarrow\infty}\bigg( 1-\underbrace{\frac{f_n^{(m-1)}}{f_n^{(m)}}(z)}_{\rightarrow\,0}\underbrace{\frac{f_n^{(m+1)}}{f_n^{(m)}}(z)}_{|\,\cdot\,|\,\leq\, M_{z}}\bigg) = 1.
		\end{align*}
		Hence, we know that there is a function $F_{m-1}\in\mathcal{M}(D\backslash\mathcal{P}_{F_m})$, such that a subsequence of $(f_n^{(m)}/f_n^{(m-1)})_n$ converges locally uniformly to $F_{m-1}$ on $D\backslash\mathcal{P}_{F_m}$.
	\end{ProofOfNo}\\~\\
	\noindent
	To show Theorem \ref{th:Schwarz}, we need an extension of Lemma \ref{le:LogarithmicDerivativesZalcman}.
	\medskip
	\begin{lemma}\label{le:Schwarz}
		Let $E\subseteq\D$ be a set with $E_\D^{(1)}=\emptyset$ and $\fn\subseteq\MD$, such that
		\begin{enumerate}
			\item[1.)] $(f_n''/f_n')_n$ converges locally uniformly to $F_1\in\mathcal{H}(\D\backslash E)$ on $\D\backslash E$.
			\item[2.)] $f_n'$ is zero-free for all $n\in\N$.
		\end{enumerate}
		Then we know that $\fn$ and $(f_n'/f_n)_n$ are quasi-normal on $\D$.\\
		If $E=\emptyset$, then $(f_n'/f_n)_n$ is normal and no subsequence of $(f_n'/f_n)_n$ ${\text{converges to }\infty}$.
	\end{lemma}
	\medskip
	\begin{proof}
		First, we want to show that $(f_n'/f_n)_n$ is quasi-normal as well and that none of its subsequences converges locally uniformly to $\infty$.\\
		Consider $(f_{n_k})_k\subseteq\fn$. Then, we can apply Lemma \ref{le:LogarithmicDerivativesZalcman} with  ${D=\D\backslash E}$ and $m=1$ to $(f_{n_k})_k$. Hence, there exists a subsequence of $(f_{n_k})_k$ (that we are going to continue to call $(f_{n_k})_k$), such that $(f_{n_k}'/f_{n_k})_k$ converges locally uniformly on $\D\backslash E$ to some $F_0\in\mathcal{M}(\D\backslash E)$. Thus, for $\mathbb{E}\coloneq E\cup\mathcal{P}_{F_0}$, we get that the sequence $(f_{n_k}'/f_{n_k})_k$ converges locally uniformly on $\D\backslash\mathbb{E}$ to the holomorphic function ${\tilde{F_0}\coloneq F_0|_{\D\backslash\mathbb{E}}}$. Now, we choose $\mathbb{E}$ as the exceptional set for our subsequence $(f_{n_k}'/f_{n_k})_k$ and get that $(f_{n_k}'/f_{n_k})_k$ does not converge to $\infty$ on $\D\backslash\mathbb{E}$. Therefore, to show our assertion on $(f_n'/f_n)_n$, we just need to show that $\mathbb{E}_\D^{(1)}=\emptyset$ or respectively $(\mathcal{P}_{F_0})_\D^{(1)}=\emptyset$.\\ \sloppy
		We assume the contrary. Then there would be a sequence $(p_j)_j\subseteq\D\backslash E$ with ${p_0\coloneq\lim_{j\rightarrow\infty}p_j\in E}$, such that $F_0(p_j)=\infty$ for all $j\in\N$ (and in particular $F_0\not\equiv0$). However, since $f_{n_k}'$ is zero-free, we know that $f_{n_k}/f_{n_k}'$ is analytic on $\D$ for all $k\in\N$ and converges uniformly to $1/\tilde{F_0}$ on some annulus within $\D$. So due to the maximum principle, $1/\tilde{F_0}$ can be analytically continued to $p_0$. However, this contradicts our assumption that $F_0(p_j)=\infty$ for all $j\in\N$.\\
		As a final step, we need to translate the convergence of $(f_n'/f_n)_n$ to $\fn$. We already know that $(f_{n_k}'/f_{n_k})_k$ is eventually bounded on $\D\backslash\mathbb{E}$. Hence, for each $w\in\D\backslash\mathbb{E}$, there is a neighborhood $U_w\subseteq\D\backslash\mathbb{E}$ and $M_w\in\R$ with
		\begin{align}\nonumber
			f_{n_k}^\#(z)=\frac{|f_{n_k}'(z)|}{1+|f_{n_k}(z)|^2}\leq \frac{2\,|f_{n_k}'(z)|}{1+|f_{n_k}(z)|}< 2\,M_w
		\end{align}
		for all $z\in U_w$ and almost all $k\in\N$.\\
		Thus, Marty's theorem implies that $\fn$ is normal on $\D\backslash\mathbb{E}$.$\hphantom{aaaaaaaaaaaaiiii}$
	\end{proof}
	\medskip\noindent
	Note that Corollary \ref{co:Extension} is very similar to Lemma \ref{le:Schwarz}. However, we still need an additional lemma to prove it.
	\medskip
	\begin{lemma}\label{le:EssentialSingularity}
		Let $D\subseteq\C$ be a domain, $z_0\in D$ and suppose that $\fn\subseteq\mathcal{M}(D)$ converges locally uniformly on $D\backslash\{z_0\}$ to ${F\in\mathcal{H}(D\backslash\{z_0\})}$. If each $f_n$ has at most $p$ poles (CM), then we know that $F$ can not have an essential singularity at $z_0$. 
	\end{lemma}
	\medskip
	\begin{proof}
		Since $F$ is analytic in $D\backslash\{z_0\}$, it follows that the set of poles of $\fn$ ${\mathcal{P}\coloneq\{z\in D;\, f_n(z)=\infty \text{ for some } n\in\N\}}$ can only accumulate in $z_0$ or in $\partial D$. However, by shrinking the domain $D$ a little bit, we are able to exclude the case of an accumulation point in $\partial D$.\\ \sloppy
		Now, we denote the poles of $f_n$ as $a_{j,n}$ (CM) (where $1\leq j\leq p$ if $f_n$ has less than $p$ poles, then we define $a_{j,n}\coloneq z_0$ for all remaining $j$). Hence, we get that  ${g_n(z)\coloneq\prod_{j=1}^p(z-a_{j,n})}$ converges locally uniformly to $(z-z_0)^p$. Thus,  ${(h_n)_n\coloneq (g_n\cdot f_n)_n\subseteq\mathcal{H}(D)}$ converges locally uniformly in $D\backslash\{z_0\}$ as well and its limit is $H(z)\coloneq(z-z_0)^p\cdot F(z)$.\\
		Next, we assume that $F$ has an essential singularity at $z_0$. Consequentially, $H$ has an essential singularity at $z_0$ too. However, this is a contradiction, since $H$ has to be holomorphic in $D$ due to the maximum principle and the Weierstraß convergence theorem.$\hphantom{mmmmmmmmmmmmmmmmmmmmmmmmmmmi}$
	\end{proof}\medskip

	\section{Proofs of the Main Results}
	\noindent
	\textit{Proof of Theorem \ref{th:Lift}:}
	\begin{ProofOfNo}
		The statement is clearly true for $j=m$.\\
		Next, we consider a subsequence $(f_{n_k})_k\subseteq\fn$ and suppose that for some ${0<j\leq m}$ there is an exceptional set $E_j$ with $(E_j)_D^{(m-j)}=\emptyset$ and a subsequence of $(f_{n_k}^{(j+1)}/f_{n_k}^{(j)})_k$ (that we are going to continue to call $(f_{n_k}^{(j+1)}/f_{n_k}^{(j)})_k$) which converges locally uniformly to some $F_j\in\mathcal{M}(D\backslash E_j)$ on $D\backslash E_j$.\\
		Now, we can apply ${ \text{Lemma } \ref{le:LogarithmicDerivativesZalcman}}$ and get that a subsequence of $(f_{n_k}^{(j)}/f_{n_k}^{(j-1)})_k$ (that we are going to continue to call $(f_{n_k}^{(j)}/f_{n_k}^{(j-1)})_k$) converges locally uniformly to some ${F_{j-1}\in\mathcal{M}\big(D\backslash(E_j\cup \mathcal{P}_{F_j})\big)}$ on $D\backslash(E_j\cup \mathcal{P}_{F_j})$. Note that $(\mathcal{P}_{F_j})_D^{(1)}\subseteq E_j$ and thus $(E_j\cup\mathcal{P}_{F_j})_D^{(m-(j-1))}=\emptyset$. Hence, we get that $(f_n^{(j)}/f_n^{(j-1)})_n$ is $Q_{m-(j-1)}$-normal.\\
		By induction, we will finally arrive at an exceptional set $E_0$ with $(E_0)_D^{(m)}=\emptyset$ and at a subsequence $(f_{n_\ell})_\ell\subseteq\fn$ where $(f_{n_\ell}'/f_{n_\ell})_\ell$ converges locally uniformly to some $F_0\in\mathcal{M}(D\backslash E_0)$ on $D\backslash E_0$. Note that by construction, we also got for $j=0,\hdots,m$ that $(f_{n_\ell}^{(j+1)}/f_{n_\ell}^{(j)})_\ell$ converges locally uniformly to some holomorphic function on $D\backslash E_0$. So the last claim of ${\text{Theorem \ref{th:Lift}}}$ follows with $U\coloneq D\backslash(P_{F_0}\cup \overline{E_0})$.
	\end{ProofOfNo}\\
	
	\noindent
	\textit{Proof of Corollary \ref{co:Extension}:}
	\begin{ProofOfNo}
		Theorem \ref{th:Lift} shows that $(f_n^{(m)}/f_n^{(m-1)})_n$ is quasi-normal and  $(f_n^{(m-1)}/f_n^{(m-2)})_n$ is $Q_2$-normal. Lemma \ref{le:LogarithmicDerivativesZalcman} has also shown it is not possible that the limit of a lifted subsequence is $\infty$. So all that is left to show is that $(f_n^{(m-1)}/f_n^{(m-2)})_n$ is quasi-normal as well.
		\sloppy
		Thus, we consider a subsequence $(f_{n_k})_{n_k}\subseteq\fn$. We know that there is a subsequence (that we are going to continue to call $(f_{n_k})_{n_k}$), such that $(f_{n_k}^{(m)}/f_{n_k}^{(m-1)})_{n_k}$ converges locally uniformly to some $F_{m-1}\in\mathcal{M}(D\backslash\mathcal{P}_{F_m})$. Similarly, we get that there is a subsequence (that we are going to continue to call $(f_{n_k})_{n_k}$), such that $(f_{n_k}^{(m-1)}/f_{n_k}^{(m-2)})_{n_k}$ converges locally uniformly on $D\backslash(\mathcal{P}_{F_m}\cup\mathcal{P}_{F_{m-1}})$. So we have to show that $(\mathcal{P}_{F_{m-1}})_D^{(1)}=\emptyset$.\\
		We assume the contrary. Then there would be a sequence $(p_j)_j\subseteq D\backslash\mathcal{P}_{F_m}$ with $\lim_{j\rightarrow\infty}p_j\eqcolon p_0\in\mathcal{P}_{F_m}$ and $F_{m-1}(p_j)=\infty$ (and in particular $F_{m-1}\not\equiv0$). We will also denote the order of the pole $p_0$ of $F_m$ as ord$(F_m,p_0)$ and consider a neighborhood $U\subseteq D$ around $p_0$, such that $\text{dist}(U,\mathcal{P}_{F_m}\backslash\{p_0\})>0$.\\
		Due to Hurwitz' theorem, we know that almost all $f_n^{(m)}$ have at most ord$(F_m,p_0)$ zeros (IM) in $U$. Therefore, we can conclude for almost all $n\in\N$ that $f_n^{(m-1)}/f_n^{(m)}$ has at most ${N\cdot\text{ord}(F_m,p_0)}$ poles (CM) in $U$. Since $(f_n^{(m-1)}/f_n^{(m)})_n$ converges locally uniformly to $1/F_{m-1}$ on $U\backslash\{p_0\}$, we can conclude that $1/F_{m-1}$ has an isolated singularity at $p_0$.\\
		Hence, we are able to apply Lemma \ref{le:EssentialSingularity} and get that $1/F_{m-1}$ has a meromorphic extension to $p_0$. Thus, $F_{m-1}$ has to be meromorphic on $U$ as well. This however contradicts our assumption that $F_{m-1}(p_n)=\infty$ for all $n\in\N$ and we get that $(\mathcal{P}_{F_{m-1}})_D^{(1)}=\emptyset$.
	\end{ProofOfNo}\\
	
	\noindent
	\textit{Proof of Theorem \ref{th:Schwarz}:}
	\begin{ProofOfNo}
		We fix $M\geq0$.\\ \sloppy
		Now suppose that $\F_M''/\F_M'\coloneq\{f''/f': f\in\F_M\}$ is not normal at $z_0\in\D$. Then, Zalcman's Lemma (Theorem \ref{ZalcmanLemma}) states that there would be $(\rho_n)_n\subseteq(0,1)$ and ${(z_n)_n\subseteq\D\backslash\{z_0\}}$ with $\lim_{n\rightarrow\infty}\rho_n=0$ and $\lim_{n\rightarrow\infty}z_n=z_0$, as well as a sequence ${(f_n''/f_n')_n\subseteq\F_M''/\F_M'}$, such that 
		\begin{align}\nonumber
			g_n(\zeta)\coloneq\frac{f_n''}{f_n'}(z_n+\rho_n\zeta		\end{align}
		converges locally uniformly to a non-constant function $g\in\mathcal{M}(\C)$.
		In particular, we get that
		\begin{align}\nonumber
			g_n'(\zeta)			=\rho_n\Bigg(\frac{f_n'''}{f_n'}-\underbrace{\bigg(\frac{f_n''}{f_n'}\bigg)^2\Bigg)(z_n+\rho_n\zeta)}_{=\,g_n^2(\zeta)}.
		\end{align}
		Since $(g_n)_n$ is eventually bounded on any compact subset of $\C\backslash\mathcal{P}_g$, we can conclude that the sequence $\Big(\mfrac{f_n'''}{f_n'}(z_n+\rho_n\,\cdot\,)\Big)_n$ converges locally uniformly to $\infty$ on ${\C\backslash(\mathcal{P}_g\cup\mathcal{Z}_{g'})}$. However, this contradicts
		\begin{align}\nonumber
			M\geq\big|S_f(z_n+\rho_n\zeta)\big| 			=\Bigg|\underbrace{\frac{f_n'''}{f_n'}(z_n+\rho_n\zeta)}_{\rightarrow\,\infty}-\frac{3}{2}\bigg(\!\underbrace{\frac{f_n''}{f_n'}(z_n+\rho_n\zeta)}_{\rightarrow\, g(\zeta)}\!\bigg)^2\Bigg|\,\,\,\,\text{ for }\zeta\in\C\backslash(\mathcal{P}_g\cup\mathcal{Z}_{g'}).
		\end{align} 
		So $\F_M''/\F_M'$ is indeed normal and any given sequence $(f_n''/f_n')_n\subseteq\mathcal\F_M''/\F_M'$ has a convergent subsequence. As a next step, we would like to show that none of these subsequences converges to $\infty$.\\
		Again, we will assume the contrary, i.e. that there is a sequence $\fn\subseteq\mathcal\F_M$, such that $(f_n''/f_n')_n$ converges locally uniformly to $\infty$ on $\D$. This would imply that
		\begin{align}\nonumber
			\lim\limits_{n\rightarrow\infty}\frac{f_n'}{f_n''}=0\quad\text{and with the Weierstraß convergence theorem}\quad \lim\limits_{n\rightarrow\infty}\frac{f_n'}{f_n''}\frac{f_n'''}{f_n''}=1.
		\end{align}
		However, this leads to
		\begin{align}\nonumber
			\lim\limits_{n\rightarrow\infty}\underbrace{\bigg(\frac{f_n'}{f_n''}\bigg)^2}_{\rightarrow\, 0}\,\,\cdot\underbrace{S_{f_n}\vphantom{\bigg(\frac{f_n'}{f_n''}\bigg)^2}}_{||\,\cdot\,||\,\leq\, M}=0
		\end{align}
		as well as
		\begin{align}\nonumber
			\lim\limits_{n\rightarrow\infty}\bigg(\frac{f_n'}{f_n''}\bigg)^2\cdot S_{f_n}			=\lim\limits_{n\rightarrow\infty}\bigg(\frac{f_n'}{f_n''}\bigg)^2\Bigg(\frac{f_n'''}{f_n'}-\frac{3}{2}\bigg(\frac{f_n''}{f_n'}\bigg)^2\Bigg)	=\underbrace{\lim\limits_{n\rightarrow\infty}\frac{f_n'}{f_n''}\frac{f_n'''}{f_n''}}_{=\,1}-\frac{3}{2}=-\frac{1}{2},
		\end{align}
		which is a contradiction.\\
		So any sequence $(h_n''/h_n')_n\subseteq\F_M''/\F_M'$ has a subsequence $(h_{n_k}''/h_{n_k}')_k$ that converges locally uniformly to some $F_1\in\MD$.\\
		Note that functions with holomorphic Schwarzians have zero-free derivatives. Thus, we can apply Lemma \ref{le:Schwarz} to $(h_{n_k})_k$ by choosing the exceptional set $E$ as the set of all poles of $F_1$ and get that $(h_n)_n$ is quasi-normal. Since  $(h_n)_n$ was chosen arbitrarily, we can conclude that $\F_M$ is quasi-normal as well.
	\end{ProofOfNo}

	\section*{Acknowledgment}
	\noindent
	I am very grateful to Jürgen Grahl for several valuable remarks.
	\bibliographystyle{amsplain}

\begin{thebibliography}{99}		
		\bibitem{Chuang} C.-T. Chuang, \textit{Normal Families of Meromorphic Functions}, World Scientific, Singapore, 1993
		\bibitem{GrahlNevo} J. Grahl, S. Nevo, Quasi-normality induced by differential inequalities, \textit{Bull. Lond. Math. Soc.} \textbf{50} (2018), 73-84
		\bibitem{Letho} O. Lehto, \textit{Univalent Functions and Teichmüller Spaces}. Springer-Verlag, New York, 1987
		\bibitem{MaMeijaMinda} W. Ma, D. Meija, D. Minda, Bounded Schwarzian and Two-Point Distortion, \textit{Comput. Methods Funct. Theory }\textbf{13} (2013), 705-715
		\bibitem{Nevo} S. Nevo, Z. Shem Tov, Differential Marty-type inequalities which lead to quasi-normality, \textit{J. Math. Anal. Appl.} \textbf{527} (2023)
		\bibitem{Zalcman} L. Zalcman, A Heuristic Principle in Complex Function Theory, \textit{Am. Math. Mon.} \textbf{82} (1975), 813-817
	\end{thebibliography}

\end{document}